\documentclass[alpha-refs]{wiley-article}

\usepackage{color}
\usepackage{ulem}

\newcommand{\eop}{\hfill $\Box$ \\ \\}
\DeclareMathOperator{\sign}{sign}

\usepackage{siunitx}

\papertype{Original Article}
\paperfield{Journal Section}

\title{Asymptotic Distribution of Constrained Nearly-Isotonic Graph Fused Lasso}

\author[1\authfn{1}]{Vladimir Pastukhov, PhD}

\affil[1]{Department of Statistics and Operations Research, University of Vienna, Austria}

\corraddress{Vladimir Pastukhov, Department of Statistics and Operations Research, University of Vienna, Austria}
\corremail{vladimir.pastukhov@univie.ac.at}

\runningauthor{Vladimir Pastukhov}

\begin{document}

\maketitle

\begin{abstract}
This paper studies the asymptotic distribution of a constrained lasso-type estimator for denoising signals defined on the nodes of a graph, where the underlying structure encodes relationships between variables. We show that, under suitable assumptions on the penalization parameters, the limiting distribution of the estimator is obtained by applying corresponding constrained procedure to the asymptotic distribution of the unrestricted estimator. Thus, the constrained estimator shares the same convergence rate the unrestricted estimator. Without the fusion penalty, the limiting distribution is obtained by applying the individual nearly isotonic estimators to the corresponding sub-vectors of the unrestricted estimator’s asymptotic distribution, similarly to the limit behaviour of isotonic regression.


\end{abstract}

\section{Introduction}
In many modern statistical applications, data are naturally structured by an underlying graph, where observations are associated with nodes and edges encode the relationships between them. A typical example is a signal with a known dependency structure, where the relationships are represented by a graph. In such settings, it is often of interest to estimate or denoise an unknown signal while incorporating structural constraints such as smoothness or monotonicity along the graph.

Let $\bm{\mathring{\beta}} \in \mathbb{R}^s$ be a vector, which may represent, for example, an unknown signal, with components indexed by $\mathcal{I} = \{\bm{i}_1, \dots, \bm{i}_s \}$. Next, suppose we have a sequence of estimators $\hat{\bm{\beta}}_{n}$ of  $\bm{\mathring{\beta}}$, for $n=1,2 \dots$, such that $\hat{\bm{\beta}}_{n}$ has the following limiting distribution
\begin{equation}\label{convbest}
n^{q}(\hat{\bm{\beta}}_{n} - \bm{\mathring{\beta}}) \stackrel{d}{\to} \bm{\psi},
\end{equation}
with $q > 0$, and where $\bm{\psi}$ is a random vector with values in $(\mathbb{R}^s, \mathcal{B}(\mathbb{R}^s))$, where $\mathcal{B}(\mathbb{R}^s)$ is the Borel $\sigma$-algebra on $\mathbb{R}^s$. In the next section, we present two illustrative examples of the base estimators. Throughout the paper, we refer to $\hat{\bm{\beta}}_{n}$ as the base estimator of $\bm{\mathring{\beta}}$, and the terms "regression", "estimator", "approximator", and "filter" are used interchangeably. We begin by introducing some notation before formulating the general problem.

Assume that we have some defined structure on the signal's support and let us start with basic definitions of partial order and isotonic regression. Let $\mathcal{I} = \{\bm{i}_{1}, \dots, \bm{i}_{s}\}$ be a set. We define the following binary relation $\preceq$ on $\mathcal{I}$.

A binary relation $\preceq$ on $\mathcal{I}$ is called partial order if 
\begin{itemize}
\item it is reflexive, i.e. $\bm{j}\preceq\bm{j}$ for all $\bm{j} \in \mathcal{I}$;
\item it is transitive, i.e. $\bm{j}_{1}, \bm{j}_{2}, \bm{j}_{3} \in \mathcal{I}$, $\bm{j}_{1} \preceq \bm{j}_{2}$ and $\bm{j}_{2} \preceq \bm{j}_{3}$ imply $\bm{j}_{1} \preceq \bm{j}_{3}$;
\item it is antisymmetric, i.e. $\bm{j}_{1}, \bm{j}_{2} \in \mathcal{I}$, $\bm{j}_{1} \preceq \bm{j}_{2}$ and $\bm{j}_{2} \preceq \bm{j}_{1}$ imply $\bm{j}_{1} = \bm{j}_{2}$.
\end{itemize}

Furthermore, a vector $\bm{\beta}\in\mathbb{R}^{s}$, indexed by $\mathcal{I}$, is called isotonic with respect to the partial order $\preceq$ on $\mathcal{I}$ if $\bm{j}_{1} \preceq \bm{j}_{2}$ implies $\beta_{\bm{j}_{1}} \leq \beta_{\bm{j}_{2}}$. We denote the set of all isotonic vectors in $\mathbb{R}^{n}$ with respect to the partial order $\preceq$ on $\mathcal{I}$ by $\bm{\mathcal{B}}^{is}$, which is a closed convex cone in $\mathbb{R}^{n}$ and it is also called isotonic cone. Next, a vector $\bm{\beta}^{I}\in \mathbb{R}^{s}$ is called the isotonic regression of an arbitrary vector $\bm{y} \in \mathbb{R}^{s}$, indexed by pre-ordered set $\mathcal{I}$, if 
\begin{eqnarray}\label{Ipo}
\bm{\beta}^{I} = \underset{\bm{\beta} \in \bm{\mathcal{B}}^{is}}{\arg \min} \sum_{\bm{j} \in \mathcal{I}}(\beta_{\bm{j}} - y_{\bm{j}})^{2}.
\end{eqnarray}

For any partial order relation $\preceq$ on $\mathcal{I}$, there exists a directed acyclic graph $G = (V,E)$, with vertex set $V = \mathcal{I}$ and the minimal set of edges $E$ such that
\begin{eqnarray*}
E = \{(\bm{j}_{1}, \bm{j}_{2}), \, \text{where} \, (\bm{j_{1}},\bm{j_{2}}) \,\text{ is the ordered pair of vertices from } \, \mathcal{I}\}.
\end{eqnarray*}
An arbitrary vector $\bm{\beta} \in \mathbb{R}^{s}$ is isotonic with respect to $\preceq$ iff $\beta_{\bm{l_{1}}} \leq \beta_{\bm{l_{2}}}$, given that $E$ contains the chain of edges from $\bm{l}_{1} \in V$ to  $\bm{l}_{2} \in V$. 

Let $D$ denote the oriented incidence matrix of the directed graph $G = (V, E)$ corresponding to the partial order $\preceq$ on $\mathcal{I}$. We choose the orientation of $D$ as follows. Suppose the graph $G$ has $n$ vertices and $m$ edges. We label the vertices by ${1, \dots, n}$ and the edges by ${1, \dots, m}$. Then $D$ is an $m \times n$ matrix defined by
\begin{equation*}\label{}
D_{i,j} =   \begin{cases}
   1, & \quad \text{if vertex $j$ is the source of edge $i$} , \\
    -1, & \quad \text{if vertex $j$ is the target of edge $i$},\\
    0, & \quad \text{otherwise}.
  \end{cases}
\end{equation*}

\begin{figure}[]
\begin{center}
\includegraphics[width=3in]{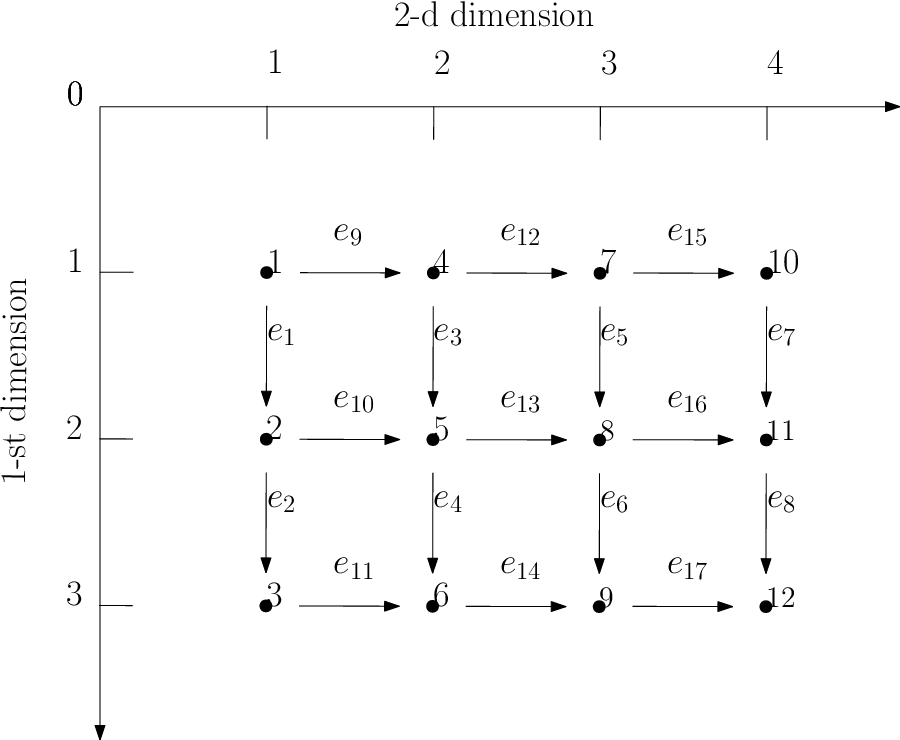}\\
\end{center}
\caption{Graph $G=(V,E)$ which corresponds to the two dimensional equally spaced grid. \label{2dftfgr}}
\end{figure}

Next, to clarify the notation, we consider a two-dimensional grid with bimonotonic constraints. The notion of bimonotonicity was first introduced in \citet{beran2010least} and is defined as follows. Let us consider the index set
\begin{eqnarray*}
\mathcal{I} = \{ \bm{i}= (i^{(1)} ,i^{(2)}): \, i^{(1)}=1,2,\dots, s_{1}, \,  i^{(2)}=1,2,\dots, s_{2}\}
\end{eqnarray*}
with the following order relation $\preceq$ on it: for $\bm{i}_{1}, \bm{i}_{2}\in \mathcal{I}$ we have $\bm{i}_{1} \preceq \bm{i}_{2}$ iff $i^{(1)}_{1} \leq i^{(1)}_{2}$ and $i^{(2)}_{1} \leq i^{(2)}_{2}$. Then, a vector $\bm{\beta}\in\mathbb{R}^{n}$, with $n=s_{1}s_{2}$, indexed by $\mathcal{I}$ is called bimonotone if it is isotonic with respect to bimonotone order $\preceq$ defined on its index $\mathcal{I}$. Further, we define the directed graph $G = (V, E)$ with vertexes $V = \mathcal{I}$, and the edges
\begin{eqnarray*}
\begin{aligned}
E ={}&  \{((l, k),(l, k+1) ): \, 1 \leq l \leq s_{1}, 1 \leq k \leq s_{2} - 1\}\\
\cup \,& \{((l, k),(l+1, k) ): \, 1 \leq l \leq s_{1}-1, 1 \leq k \leq s_{2} \}.
\end{aligned}
\end{eqnarray*}
The oriented graph for $3\times 4$ grid is displayed on Figure \ref{2dftfgr}. The oriented incidence matrix $D \in \mathbb{R}^{17\times 12}$ corresponding to the graph $G=(V,E)$, displayed in Figure (\ref{2dftfgr}), is given in (\ref{2dDmat}).

Now we can define isotonic regression equivalently to (\ref{Ipo}). The isotonic regression of an arbitrary vector $\bm{y} \in \mathbb{R}^{s}$,  indexed by the partially ordered index set $\mathcal{I}$, is given by
\begin{eqnarray*}\label{}
\bm{\beta}^{I} = \underset{\bm{\beta}}{\arg \min} \sum_{\bm{j} \in \mathcal{I}}(\beta_{\bm{j}} - y_{\bm{j}})^{2},
\end{eqnarray*}
subject to $\beta_{\bm{l_{1}}} \leq \beta_{\bm{l_{2}}}$, whenever $E$ contains the chain of the edges from $\bm{l}_{1} \in V$ to  $\bm{l}_{2} \in V$.  

For a general introduction to the subject of constrained inference we refer to the monographs: \cite{barlowstatistical},  \cite{groeneboom}, \cite{robertson1988order},  and \cite{silvapsen}. These monographs consider isotonic regression in various settings, addressing fundamental questions such as existence, uniqueness and computational aspects of the estimator. The limiting distribution of general isotonic regression has been studied in the discrete and continuous settings in \cite{anevski2018asymptotic}, and \cite{han2020limit}, respectively.

\begin{equation}\label{2dDmat}
D = \begin{pmatrix}
1 & -1 & 0 & 0 & 0 &0 & 0 & 0 & 0& 0 & 0& 0\\
0 & 1 & -1 & 0 & 0 &0 & 0 & 0 & 0& 0 & 0& 0\\
0 & 0 & 0 & 1 & -1 & 0 & 0 & 0 & 0& 0 & 0& 0\\
0 & 0 & 0 & 0 & 1 &-1 & 0 & 0 & 0& 0 & 0& 0\\
0 & 0 & 0 & 0 & 0 & 0 & 1 & -1 & 0& 0 & 0& 0\\
0 & 0 & 0 & 0 & 0 & 0 & 0 & 1 & -1& 0& 0& 0\\
0 & 0 & 0 & 0 & 0 & 0 & 0 & 0 & 0 & 1 & -1& 0\\
0 & 0 & 0 &0 & 0 & 0 & 0 & 0 & 0 & 0 & 1& -1\\
1 & 0 & 0 & -1 & 0 & 0 & 0 & 0 & 0 & 0& 0& 0\\
0 & 1 & 0 & 0 & -1 & 0 & 0 & 0 & 0 & 0& 0& 0\\
0 & 0 & 1 & 0 & 0 & -1 & 0 & 0 & 0 & 0& 0& 0\\
0 & 0 & 0 & 1 & 0 & 0 & -1 & 0 & 0 & 0& 0& 0\\
0 & 0 & 0 & 0 & 1 & 0 & 0 & -1 & 0 & 0& 0& 0\\
0 & 0 & 0 & 0 & 0 & 1 & 0 & 0 & -1 & 0& 0& 0\\
0 & 0 & 0 & 0 & 0 & 0 & 1 & 0 & 0 & -1& 0& 0\\
0 & 0 & 0 & 0 & 0 & 0 & 0 & 1 & 0 &0& 0& -1
\end{pmatrix}
\end{equation}

A relaxed form of isotonic regression, called nearly isotonic regression, allows violations of the monotonicity constraints and is defined by
\begin{equation*}\label{}
\hat{\bm{\beta}}^{NI}(\bm{y}, \lambda^{NI}) = \underset{\bm{\beta} \in \mathbb{R}^{n}}{\arg \min} \, \frac{1}{2} ||\bm{y} - \bm{\beta}||_{2}^{2} + \lambda^{NI}\sum_{(\bm{i},\bm{j})\in E}|\beta_{\bm{i}} - \beta_{\bm{j}}|_{+},
\end{equation*}
where $x_{+} =  x \cdot 1\{x > 0 \}$. Alternatively, with the use of the incidence matrix $D$, it can be written as 
\begin{equation*}\label{}
\begin{aligned}
\hat{\bm{\beta}}^{NI}(\bm{y}, \lambda^{NI}) &= {} 
 \underset{\bm{\beta} \in \mathbb{R}^{n}}{\arg \min} \, \frac{1}{2} ||\bm{y} - \bm{\beta}||_{2}^{2} +  \lambda^{NI}||D\bm{\beta}||_{+},
\end{aligned}
\end{equation*}
where $||D\bm{x}||_{+} = \sum_{k=1}^{s} (x_{k})_{+}$, for $\bm{x} \in \mathbb{R}^{s}$.

Nearly-isotonic regression was first introduced by \cite{tibshirani2011nearly}. Its one-dimensional limiting distribution was established in \cite{minami2020estimating}, and properties for general graphs, including computational aspects, were studied in \citet{pastukhov2022fused}.

Next, if the goal is the overall denoising without order control, we use graph fused lasso estimator, which is given by
\begin{equation*}\label{}
\hat{\bm{\beta}}^{F}(\bm{y}, \lambda^{F}) = \underset{\bm{\beta} \in \mathbb{R}^{n}}{\arg \min} \, \frac{1}{2} ||\bm{y} - \bm{\beta}||_{2}^{2} + \lambda^{F}\sum_{(\bm{i},\bm{j})\in E}|\beta_{\bm{i}} - \beta_{\bm{j}}|,
\end{equation*}
or, equivalently, we can write it as
\begin{equation*}\label{}
\begin{aligned}
\hat{\bm{\beta}}^{F}(\bm{y}, \lambda^{F}) &= {} 
 \underset{\bm{\beta} \in \mathbb{R}^{n}}{\arg \min} \, \frac{1}{2} ||\bm{y} - \bm{\beta}||_{2}^{2} +  \lambda^{F}||D\bm{\beta}||_{1}.
\end{aligned}
\end{equation*}
Fused lasso on a grid was first introduced in signal processing by \cite{rudin1992nonlinear}. The solution path for a general graph was obtained in \cite{hoefling2010path}. An exact solution for lattice data with different penalization parameters, based on the taut-string algorithm, was proposed in \cite{barbero2018modular}. Recent results on fused lasso estimation over graphs include \cite{chen2023more, padilla2017dfs, padilla2022variance}. The asymptotic properties in the one-dimensional and general settings were studied in \citet{rinaldo2009properties} and \cite{tibshirani2005sparsity}, respectively. 

\section{Statement of the problem}
We propose the following improved estimator of $\bm{\mathring{\beta}}$, based on a given base estimator $\hat{\bm{\beta}}_{n}$ that satisfies (\ref{convbest}) in the Introduction:

\begin{equation}\label{FNIGc}
\begin{aligned}
\hat{\bm{\beta}}^{*}_{n} 
&=&& \underset{\bm{\beta} \in \mathbb{R}^{n}}{\arg\min}  
\frac{1}{2} \big\| \hat{\bm{\beta}}_{n} - \bm{\beta} \big\|_{2}^{2} 
+ \lambda_{n}^{F} \sum_{(\bm{i}, \bm{j}) \in E} 
\big| \beta_{\bm{i}} - \beta_{\bm{j}} \big|
+ \lambda_{n}^{NI} \sum_{(\bm{i}, \bm{j}) \in E} 
\big| \beta_{\bm{i}} - \beta_{\bm{j}} \big|_{+}
, \\[0.5em]
\text{s.t.} 
& &&\sum_{i=1}^{s} \hat{\beta}_{n,i} = \sum_{i=1}^{s} \beta_{i}, 
\quad 
\min(\hat{\bm{\beta}}_{n})\, \bm{1} 
\leq \bm{\beta}
\leq \max(\hat{\bm{\beta}}_{n})\, \bm{1}.
\end{aligned}
\end{equation}

We extend the approach proposed by \citet{pastukhov2022fused} and investigate the asymptotic properties of the resulting estimator. To clarify the problem formulation and highlight the role and importance of the additional constraints, we provide the following illustrative examples, which are the standard examples in the subject of constraint inference cf. \cite{robertson1988order}.

In the first example, we consider the problem of estimation a bimonotone signal in the presence of noise. Suppose that we observe a signal on a grid $Y_{\bm{i}}$ indexed by $\mathcal{I} = \{ \bm{i}= (i^{(1)} ,i^{(2)}): \, i^{(1)}=1,2,\dots, s_{1}, \,  i^{(2)}=1,2,\dots, s_{2}\}$ with bimonotone order $\preceq$ defined on $\mathcal{I}$. Next, assume that $\bm{i}_k \in \mathcal{I}$, for $ k=1, \dots, n$, and that we observe  data of the form
\begin{eqnarray*}
Y_{\bm{i}_k} = \mathring{g}_{\bm{i}_k} + \varepsilon_{\bm{i}_k}, \quad k=1, \dots, n,
\end{eqnarray*}
where $\mathring{g} \in \mathbb{R}^{s_1\times s_2}$ denotes the underlying signal, $\varepsilon_{\bm{i}_k}$ is a sequence of  i.i.d. random variables with $\mathbb{E}[\varepsilon_{\bm{i}_k}] = 0$, $\mathrm{Var}[\varepsilon_{\bm{i}_k}] = \sigma^{2} < \infty$.  
Further, for $\bm{j} \in \mathcal{I}$ let 
\begin{eqnarray*}\label{}
\hat{g}_{n, \bm{j}} =  \frac{\sum_{k=1}^{n} Y_{\bm{i}_k} 1 \{ \bm{i}_k = \bm{j} \}} {\sum_{k=1}^{n} 1 \{ \bm{i}_k = \bm{j} \}},
\end{eqnarray*}
and let 
\begin{eqnarray}\label{wxn}
w_{\bm{x}}^{(\bm{j})} = \frac{\sum_{i=1}^{n} 1\{\bm{i}_k = \bm{j}\}}{n},
\end{eqnarray}
assuming 
\begin{eqnarray*}\label{}
\bm{w}^{(n)} \to \bm{w},
\end{eqnarray*}
as $n \to \infty$, where $\bm{w}^{(n)}$ is a sequence of vectors in $\mathbb{R}_{+}^{s_{1} \times s_{2}}$ whose components are from (\ref{wxn}), and $\bm{w} \in \mathbb{R}_{+}^{s_{1} \times s_{2}}$. Then, we have
\begin{eqnarray}\label{asunreg}
n^{1/2}(\hat{\bm{g}}_{n} - \bm{\mathring{g}}) \stackrel{d}{\to} \bm{Y}_{\bm{0}, \Sigma},
\end{eqnarray}
where $\bm{Y}_{\bm{0}, \Sigma}$ is a Gausian vector with mean zero and diagonal covariance matrix $\Sigma$, whose elements are given by $\Sigma _{kk} = \sigma^{2} w_{k}$, for $k= 1, \dots, s_{1} \times s_{2}$. This example can be viewed as a simplified instance of a general nonlinear least squares estimator, which has been studied, for example, in \cite{jennrich1969asymptotic}, and \cite{wu1981asymptotic}.

Furthermore, the estimator $\hat{\bm{g}}_{n}$ of $\bm{\mathring{g}}$ can be smoothed (or denoised) using a fused lasso penalty. If, in addition, the underlying signal is assumed to be bimonotone, it can be further regularised using a nearly-isotonic term to enforce approximate monotonicity. Finally, equality and inequality constraints can be imposed as defined in (\ref{FNIGc}). 

In general, in physics applications, one may need to smooth an intensity plot. In this setting, additional constraints can ensure that the final estimator preserves the total energy and that the filtered signal remains within the range of the base estimator.

In the second example, we consider the problem of estimation of probability mass function with the support with bimonotone structure . Suppose that we have observed $Z_{1}, Z_{2}, \dots, Z_{n}$ i.i.d. random variables taking values in $\mathcal{I} = \{ \bm{i}= (i^{(1)} ,i^{(2)}): \, i^{(1)}=1,2,\dots, s_{1}, \,  i^{(2)}=1,2,\dots, s_{2}\}$ with probability mass function $\bm{p} \in \mathbb{R}^{s_{1}\times s_{2}}$. The empirical estimator of $\bm{p}$ is then given by
the empirical estimator of $\bm{p}$ is then given by
\begin{eqnarray*}\label{}
	\hat{p}_{n, \bm{i}} = \frac{n_{\bm{i}}}{n}, \quad n_{\bm{i}}= \sum_{j=1}^{n} 1 \{ Z_{j} = \bm{i} \}, \quad \text{$ \bm{i} \in \mathcal{I}$},
\end{eqnarray*}
The empirical estimator is asymptotically Gaussian
\begin{eqnarray*}\label{}
n^{1/2}(\hat{\bm{p}}_{n} - \bm{p}) \stackrel{d}{\to} \bm{Y}_{\bm{0}, C},
\end{eqnarray*}
where $\bm{Y}_{\bm{0}, C}$ is a Gaussian vector in $\mathbb{R}^{s_{1}\times s_{2}}$, with mean zero and the covariance matrix $C  = p_{\bm{i}} \delta_{\bm{i},\bm{i'}} - p_{\bm{i}}p_{\bm{i'}}$,  where $\delta_{\bm{i}\bm{j}}=1$, if $\bm{i}=\bm{j}$ and $0$ otherwise.

Next, we can apply the procedure (\ref{FNIGc}) to smooth the empirical estimator $\hat{\bm{p}}_{n}$, which serves as the base estimator $\hat{\bm{\beta}}_{n}$, and $\bm{\psi}$ in (\ref{convbest}) is equal to $\bm{Y}_{\bm{0}, C}$. Therefore, it is natural in this setting to impose both equality constraints and non-negativity constraints on the final smoothed estimator $\hat{\bm{\beta}}^{*}_{n}$. With a slight abuse of notation, we treat the empirical estimator as the “signal” defined on a graph-structured support.
The isotonisation of the empirical estimator was studied in detail in \citet{jankowski2009estimation}. A potential application of the estimator is the monotonic smoothing of a bivariate histogram of population-based body mass index (BMI) deciles for mothers and fathers of a group of children with obesity \cite{hebebrand2000epidemic}.

In both examples above, we have $q = 2$. In this paper, we consider the more general case $q > 0$, which is motivated by \cite{andrews1999estimation}, where constrained parametric estimators were studied and examples with $q \neq 2$ were provided. Moreover, in \cite{anevski2018asymptotic}, the authors studied isotonic regression estimators over a general countable pre-ordered set under general assumptions on $q$.

\section{Asymptotic distribution of the estimator}
We start with the proof of a property shared by both fused and nearly-isotonic estimators.
\begin{theorem}\label{prop_constr}
For any values of the penalization parameters $\lambda^{F}_{n}$ and $\lambda^{NI}_{n}$, the solution to the optimization problem in (\ref{FNIGc}) is equal to the solution of the following problem
\begin{equation}\label{FNIGc1}
\begin{aligned}
\hat{\bm{\beta}}^{*}_{n} 
= \underset{\bm{\beta} \in \mathbb{R}^{n}}{\arg\min}  
\frac{1}{2} \big\| \hat{\bm{\beta}}_{n} - \bm{\beta} \big\|_{2}^{2} 
+ \lambda_{n}^{F} \sum_{(\bm{i}, \bm{j}) \in E} 
\big| \beta_{\bm{i}} - \beta_{\bm{j}} \big|
+ \lambda_{n}^{NI} \sum_{(\bm{i}, \bm{j}) \in E} 
\big| \beta_{\bm{i}} - \beta_{\bm{j}} \big|_{+},
\end{aligned}
\end{equation}
i.e.,  the additional constraints 
\begin{equation}\label{FNIGcond}
\begin{aligned}
\sum_{i=1}^{s} \hat{\beta}_{n,i} = \sum_{i=1}^{s} \beta_{i}, 
\quad 
\min(\hat{\bm{\beta}}_{n})\, \bm{1} 
\leq \bm{\beta} 
\leq \max(\hat{\bm{\beta}}_{n})\, \bm{1},
\end{aligned}
\end{equation}
are redundant and automatically satisfied.
\end{theorem}

\textbf{Proof of Theorem \ref{prop_constr}} We will show that the solution to the problem (\ref{FNIGc1}) satisfies to the conditions (\ref{FNIGcond}). First, for simplicity of notation, we suppress the subscript $n$ in $\hat{\bm{\beta}}_{n}$, $\hat{\bm{\beta}}^{*}_{n}$, $\lambda^{F}_{n}$, and in $\lambda^{NI}_{n}$ Next, analogously to the approach in Section 2 of \citet{hoefling2010path}, assume that, for some penalization parameters $\lambda^{F}$ and $\lambda^{NI}$, we have $k_{F}$ fused constant regions $F_{1}, \dots, F_{k_{F}}$ in the graph $G = (V, E)$, i.e. the subsets $F_{1}, \dots, F_{k_{F}}$ of the underlying index set $\mathcal{I}$, such that
\begin{itemize}
\item $\cup_{i=1}^{n_{F}} F_{i} = \mathcal{I}$,
\item $F_{i}\cap F_{j} = \emptyset$, for $i\neq j$,
\item if $\bm{u}, \bm{v} \in F_{i}$, then  $\bm{u}$ and  $\bm{v}$ are connected in $G$ by only going over the nodes in the subgraph of $G$ induced by $F_{i}$,
\item if $\bm{u}, \bm{v} \in F_{i}$, then $\hat{\beta}^{*}_{\bm{u}} = \hat{\beta}^{*}_{\bm{v}}$, and if $\bm{u} \in F_{i}$, $\bm{v} \in F_{j}$, $i \neq j$, and $F_{i}$ and $F_{j}$ have a connecting edge, then $\hat{\beta}^{*}_{\bm{u}} \neq  \hat{\beta}^{*}_{\bm{v}}$.
 \end{itemize}
 
We set $\beta_{\bm{l}} = \beta_{F_{i}}$ for all $\bm{l} \in F_{i}$, and let $\big| \{ (\bm{u},\bm{v})\in E: \bm{u}\in F_{i}, \bm{v}\in F_{j} \} \big|$ denote the cardinality of a set of pairs  $\{ (\bm{u},\bm{v})\in E: \bm{u}\in F_{i}, \bm{v}\in F_{j} \}$. Then, the problem (\ref{FNIGc1}) can be written as 
\begin{equation*}\label{}
\begin{aligned}
\hat{\bm{\beta}}^{*}
&= \underset{\bm{\beta} \in \mathbb{R}^{k_{F}}}{\arg\min}  
\frac{1}{2}\sum_{i=1}^{k_{F}}\sum_{\bm{l} \in F_{i}} ( \hat{\beta}_{\bm{l}} - \beta_{F_{i}} )_{2}^{2} 
+ \lambda^{F} \sum_{i < j}  \big|  \{ (\bm{u},\bm{v})\in E: \bm{u}\in F_{i}, \bm{v}\in F_{j} \}\big|
| \beta_{F_{i}} - \beta_{F_{j}} | \\
&+ \lambda^{NI} \sum_{i < j} \big| \{ (\bm{u},\bm{v})\in E: \bm{u}\in F_{i}, \bm{v}\in F_{j} \} \big|
| \beta_{F_{i}} - \beta_{F_{j}} |_{+}.
\end{aligned}
\end{equation*}
Next, the derivative of objective function with respect to $\beta_{F_{i}}$ at $\hat{\beta}^{*}_{F_{i}}$, for $i =1, \dots, k_{F}$ is $0$, i.e.  we have
\begin{equation*}\label{}
\begin{aligned}
|F_{i}|\hat{\beta}^{*}_{F_{i}} - \sum_{j \in F_{i}} \hat{\beta}_{j} &+ \lambda^{F} \sum_{i < j} \big| \{ (\bm{u},\bm{v})\in E: \bm{u}\in F_{i}, \bm{v}\in F_{j} \} \big|
\sign( \hat{\beta}^{*}_{F_{i}} - \hat{\beta}^{*}_{F_{j}} )\\ 
&+  \lambda^{NI} \sum_{i < j} \big| \{ (\bm{u},\bm{v})\in E: \bm{u}\in F_{i}, \bm{v}\in F_{j}\} \big|
1\{ \hat{\beta}^{*}_{F_{i}} - \hat{\beta}^{*}_{F_{j}} > 0\} = 0,
\end{aligned}
\end{equation*}
and, consequently,
\begin{equation*}\label{}
\begin{aligned}
\hat{\beta}^{*}_{F_{i}} = \frac{\sum_{j \in F_{i}} \hat{\beta}_{j}}{|F_{i}|} &- \frac{\lambda^{F}}{|F_{i}|} \sum_{i < j} \big| \{  (\bm{u},\bm{v})\in E: \bm{u}\in F_{i}, \bm{v}\in F_{j} \} \big|
\sign( \hat{\beta}^{*}_{F_{i}} - \hat{\beta}^{*}_{F_{j}} )\\ 
&-  \frac{\lambda^{NI}}{|F_{i}|} \sum_{i < j} \big| \{  (\bm{u},\bm{v})\in E: \bm{u}\in F_{i}, \bm{v}\in F_{j} \} \big|
1\{ \hat{\beta}^{*}_{F_{i}} - \hat{\beta}^{*}_{F_{j}} > 0\}. 
\end{aligned}
\end{equation*}

Further, assume that $\hat{\beta}^{*}_{F_{i}}$ is the smallest. Then, we have
\begin{equation*}\label{}
\begin{aligned}
&\big|\{(\bm{u},\bm{v})\in E: \bm{u}\in F_{i}, \bm{v}\in F_{j}\}\big|
\sign( \hat{\beta}_{F_{i}} - \hat{\beta}_{F_{j}} ) < 0,\\ 
&\big|\{ (\bm{u},\bm{v})\in E: \bm{u}\in F_{i}, \bm{v}\in F_{j}\}
\big| 1\{ \hat{\beta}_{F_{i}} - \hat{\beta}_{F_{j}} > 0\} = 0,
\end{aligned}
\end{equation*}
and, therefore, 
\begin{equation*}\label{}
\begin{aligned}
\hat{\beta}_{F_{i}} \geq  \frac{\sum_{j \in F_{i}}  \hat{\beta}_{j} }{|F_{i}|} \geq \min(\hat{\bm{\beta}}).
\end{aligned}
\end{equation*}
The proof for the upper bound can be done in the same way. Assume that $\hat{\beta}^{*}_{F_{t}}$ is the largest. Then, in this case we have 
\begin{equation*}\label{}
\begin{aligned}
&\big|\{(\bm{u},\bm{v})\in E: \bm{u}\in F_{t}, \bm{v}\in F_{j}\} \big|
\sign( \hat{\beta}_{F_{t}} - \hat{\beta}_{F_{j}} ) > 0,\\ 
&\big| \{ (\bm{u},\bm{v})\in E: \bm{u}\in F_{t}, \bm{v}\in F_{j}\} \big|
1\{ \hat{\beta}_{F_{t}} - \hat{\beta}_{F_{j}} > 0\} > 0,
\end{aligned}
\end{equation*}
and, consequently, 
\begin{equation*}\label{}
\begin{aligned}
\hat{\beta}_{F_{t}} \leq  \frac{\sum_{j \in F_{i}}  \hat{\beta}_{j} }{|F_{i}|} \leq \max(\hat{\bm{\beta}}).
\end{aligned}
\end{equation*}
Next, the equality constraint follows from Theorem 2..4 in \citet{pastukhov2025fused}, and, finally the result of the theorem follows.
\eop

We first note that the range condition in  (\ref{FNIGcond})  is an expected property of the fused estimator. Isotonic regression also satisfies both additional constraints (cf. Theorems 1.3.3 and 1.3.4 in \citet{robertson1988order}, respectively), but the proofs rely on the fact that the isotonic regression estimator is the projection onto a closed convex cone, which is not the case in the problem we consider in this paper. 
\begin{corollary}\label{}
Assume that the base estimator $\hat{\bm{\beta}}_{n}$ is a probability vector (for example, the empirical estimator), then, the solution to (\ref{FNIGc1}) is also a probability vector.
\end{corollary}
In the next Theorem we derive the asymptotic distribution of fused nearly-isotonic regression. 
\begin{theorem}\label{sol_fnia}
Assume that $\lambda^{F}_{n}/n^{q}\to \lambda^{F}_{0}$ and $\lambda^{NI}_{n}/n^{q}\to \lambda^{NI}_{0}$ then
\begin{equation*}\label{}
n^{q}(\hat{\bm{\beta}}^{*}_{n} - \bm{\mathring{\beta}}) \stackrel{d}{\to} \arg \min V(\bm{w}),
\end{equation*}
where $\bm{w} \in \mathbb{R}^{s}$, and
\begin{equation*} 
\begin{aligned}
V(\bm{w}) &= -2\bm{\psi}^{T}\bm{w} + \bm{w}^{T}\bm{w} &\\
&+\lambda^{F}_{0} \sum_{(\bm{i}, \bm{j}) \in E} (w_{\bm{i}} -w_{\bm{j}})
\sign( \mathring{\beta}_{\bm{i}} - \mathring{\beta}_{\bm{j}})1\{ \mathring{\beta}_{\bm{i}} \neq \mathring{\beta}_{\bm{j}}\} + \lambda^{F}_{0} \sum_{(\bm{i}, \bm{j}) \in E} |w_{\bm{i}} -w_{\bm{j}}|
1\{ \mathring{\beta}_{\bm{i}} = \mathring{\beta}_{\bm{j}}\}\\
&+\lambda^{NI}_{0} \sum_{(\bm{i}, \bm{j}) \in E} (w_{\bm{i}} -w_{\bm{j}})
1\{ \mathring{\beta}_{\bm{i}} > \mathring{\beta}_{\bm{j}}\} + \lambda^{NI}_{0} \sum_{(\bm{i}, \bm{j}) \in E} |w_{\bm{i}} -w_{\bm{j}}|_{+}
1\{ \mathring{\beta}_{\bm{i}} = \mathring{\beta}_{\bm{j}}\},
\end{aligned}
\end{equation*}
and $\bm{\psi}$ is the limiting distribution of the base estimator $\hat{\bm{\beta}}_{n}$.
\end{theorem}

\textbf{Proof of Theorem \ref{sol_fnia}} First, let us define $V_{n}(\bm{w})$ by
\begin{equation*} 
\begin{aligned}
V_{n}(\bm{w}) &= ||n^{q}(\bm{\hat{\beta}}_{n} - \bm{\mathring{\beta}})  - \bm{w}||_{2}^{2} - ||n^{q}(\bm{\hat{\beta}}_{n} - \bm{\mathring{\beta}})||_{2}^{2} \\
&+ \lambda^{F}_{n}\Big\{\sum_{(\bm{i}, \bm{j}) \in E} | \mathring{\beta}_{\bm{i}} - \mathring{\beta}_{\bm{j}} + \frac{w_{\bm{i}} -w_{\bm{j}}}{n^{q}}| -  | \mathring{\beta}_{\bm{i}} - \mathring{\beta}_{\bm{j}}|\Big\} + \lambda^{NI}_{n}\Big\{\sum_{(\bm{i}, \bm{j}) \in E} | \mathring{\beta}_{\bm{i}} - \mathring{\beta}_{\bm{j}} + \frac{w_{\bm{i}} -w_{\bm{j}}}{n^{q}}|_{+} -  | \mathring{\beta}_{\bm{i}} - \mathring{\beta}_{\bm{j}}|_{+} \Big\},
\end{aligned}
\end{equation*}
where $\bm{w} \in \mathbb{R}^{s}$. 

Second, $V_{n}(\bm{w})$ is minimazed at $\hat{\bm{w}}_{n} = n^{q}(\hat{\bm{\beta}}_{n}^{*} - \bm{\mathring{\beta}})$, and 
\begin{equation*} 
\begin{aligned}
V_{n}(\bm{w}) &= ||n^{q}(\bm{\hat{\beta}}_{n} - \bm{\mathring{\beta}})  - \bm{w}||_{2}^{2} - ||n^{q}(\bm{\hat{\beta}}_{n} - \bm{\mathring{\beta}})||_{2}^{2} \stackrel{d}{\to} -2\bm{\lambda}^{T}\bm{w} + \bm{w}^{T}\bm{w}.
\end{aligned}
\end{equation*}

Third, we have 
\begin{equation*} 
\begin{aligned}
 &\lambda^{F}_{n}\Big\{\sum_{(\bm{i}, \bm{j}) \in E} | \mathring{\beta}_{\bm{i}} - \mathring{\beta}_{\bm{j}} + \frac{w_{\bm{i}} -w_{\bm{j}}}{n^{q}}| -  | \mathring{\beta}_{\bm{i}} - \mathring{\beta}_{\bm{j}}|\Big\} \\
 &\to \lambda^{F}_{0} \sum_{(\bm{i}, \bm{j}) \in E} (w_{\bm{i}} -w_{\bm{j}})
\sign( \mathring{\beta}_{\bm{i}} - \mathring{\beta}_{\bm{j}})1\{ \mathring{\beta}_{\bm{i}} \neq \mathring{\beta}_{\bm{j}}\} + \lambda^{F}_{0} \sum_{(\bm{i}, \bm{j}) \in E} |w_{\bm{i}} -w_{\bm{j}}|
1\{ \mathring{\beta}_{\bm{i}} = \mathring{\beta}_{\bm{j}}\},
\end{aligned}
\end{equation*}
and
\begin{equation*} 
\begin{aligned}
& \lambda^{NI}_{n}\Big\{\sum_{(\bm{i}, \bm{j}) \in E} | \mathring{\beta}_{\bm{i}} - \mathring{\beta}_{\bm{j}} + \frac{w_{\bm{i}} -w_{\bm{j}}}{n^{q}}|_{+} -  | \mathring{\beta}_{\bm{i}} - \mathring{\beta}_{\bm{j}}|_{+} \Big\}\\
&\to \lambda^{NI}_{0} \sum_{(\bm{i}, \bm{j}) \in E} (w_{\bm{i}} -w_{\bm{j}})
1\{ \mathring{\beta}_{\bm{i}} > \mathring{\beta}_{\bm{j}}\} + \lambda^{NI}_{0} \sum_{(\bm{i}, \bm{j}) \in E} |w_{\bm{i}} -w_{\bm{j}}|_{+}
1\{ \mathring{\beta}_{\bm{i}} = \mathring{\beta}_{\bm{j}}\}.
\end{aligned}
\end{equation*}

Therefore, we have shown that 
\begin{equation*} 
\begin{aligned}
V_{n}(\bm{w}) \stackrel{d}{\to} V(\bm{w}).
\end{aligned}
\end{equation*}
Next, $V_{n}(\bm{w})$ is convex and $V(\bm{w})$ has unique minimum. The result of the Theorem now follows from Theorem 3.2 in \citet{geyer1996asymptotics}. This result is, in some sense, analogous to the asymptotic results for lasso in regression setup, cf. \citet{knight2000asymptotics, tibshirani2005sparsity}. \eop

The asymptotic behaviour of the estimator introduced in this paper is similar to the asymptotic behaviour of isotonic regression. Indeed, as shown in \citet{jankowski2009estimation}, the asymptotic distribution of the isotonic regression is given by the concatenation of the separate isotonic regressions of the certain subvectors of an unrestrecred estimator’s asymptotic distribution. 
 
Let us assume that $\lambda^{F}_{n} = 0$, $\lambda^{NI}_{n}/n^{q}\to \lambda^{NI}_{0}$, and that the true signal is isotonic. In this case, the asymptotic distribution is given by
\begin{equation*}
n^{q}(\hat{\bm{\beta}}^{*}_{n} - \bm{\mathring{\beta}}) \stackrel{d}{\to} \arg \min V(\bm{w}),
\end{equation*}
where
\begin{equation*} 
\begin{aligned}
V(\bm{w}) &= -2\bm{\psi}^{T}\bm{w} + \bm{w}^{T}\bm{w} + \lambda^{NI}_{0} \sum_{(\bm{i}, \bm{j}) \in E} |w_{\bm{i}} -w_{\bm{j}}|_{+}
1\{ \mathring{\beta}_{\bm{i}} = \mathring{\beta}_{\bm{j}}\}\\ 
\end{aligned}
\end{equation*}
or, equivalently, we have 
\begin{equation*}
n^{q}(\hat{\bm{\beta}}^{*}_{n} - \bm{\mathring{\beta}}) \stackrel{d}{\to} \arg \min V^{'}(\bm{w}),
\end{equation*}
where 
\begin{equation*} 
\begin{aligned}
V^{'}(\bm{w}) &=  ||\bm{w} - \bm{\psi}||_{2}^{2} + \lambda^{NI}_{0} \sum_{(\bm{i}, \bm{j}) \in E} |w_{\bm{i}} -w_{\bm{j}}|_{+},
1\{ \mathring{\beta}_{\bm{i}} = \mathring{\beta}_{\bm{j}}\},\\ 
\end{aligned}
\end{equation*}
which, in words, means a concatenation of separate nearly-isotonic regressions over constant regions of the true signal, applied to the asymptotic distribution $\bm{\psi}$ of the base estimator. 
 
 \section{Conclusions}
In this paper, we studied asymptotic properties of constrained fused lasso nearly-isotonic regression. We derived the asymptotic distribution for the case of a general directed acyclic graph and found that the asymptotic behaviour of nearly-isotonic regression is similar to the behaviour of isotonic regression.

One of the possible directions is to study asymptotic properties of nearly-isotonic regression with, for example, $\ell_{p}$ loss, with $p\neq 2$. Next, to our knowledge, the asymptotic distribution of lasso-type estimators with data-driven penalization parameters remains an open problem. Furthermore, it is important to study the behaviour of the penalised estimator when the size of the underlying signal is not fixed but depends on the sample size.

\section*{conflict of interest}
The author declares that there is no conflict of interest.

\bibliography{sample.bib}

\end{document}